%% file: PlaneOProots.tex
\pdfoutput=1
\documentclass{article}
\usepackage{amssymb,amsmath}
\usepackage{wrapfig}
\usepackage{mathrsfs}
\usepackage{amsbsy}
\usepackage{color}
\usepackage{graphicx}

\def\wt{\widetilde}
\def\wh#1{\mathrm{Co}\!\left(#1\right)}
\def\ds{\displaystyle}
\def\supp{\mathrm {supp}}

\makeatletter
\@addtoreset{equation}{section}
\makeatother
\def\tr{\mathrm {Tr}}
\def\le{\left}
\def\ri{\right}
\def\QED{{\bf Q.E.D.}\par\vskip 5pt}
\def\bc{\begin{corollary}}
\def\ec{\end{corollary}}
\def\&{&{\hskip -20pt}}
\def\ov{\overline}
\def\br{\begin{remark}\rm\small}
\def\1{{\bf 1}}
\def\er{\end{remark}}
\def\bt{\begin{theorem}}
\def\et{\end{theorem}}

\def\bx{\begin{example}}
\def\ex{\end{example}}
\def\bd{\begin{definition}}
\def\ed{\end{definition}}
\def\bp{\begin{proposition}\rm}
\def\bl{\begin{lemma}\em}
\def\el{\end{lemma}}
\def\ep{\end{proposition}}
\def\bea{\begin{eqnarray}}
\def\eea{\end{eqnarray}}
\def \pa{\partial}
\def\C{{\mathbb C}}
\def\R{{\mathbb R}}
\def\N{{\mathbb N}}

\newtheorem{theorem}{Theorem}[section]
\newtheorem{example}{Example}[section]
\newtheorem{coroll}{Corollary}[section]

\newtheorem{lemma}{Lemma}[section]
\newtheorem{remark}{Remark}[section]

\newtheorem{proposition}{Proposition}[section] 
\newtheorem{definition}{Definition}[section]
\def\br{\begin{remark}}
\def\er{\end{remark}}
\def\bt{\begin{theorem}}
\def\et{\end{theorem}}
\def\bc{\begin{coroll}}
\def\ec{\end{coroll}}
\def\bl{\begin{lemma}}
\def\el{\end{lemma}}
\def\bd{\begin{definition}}
\def\ed{\end{definition}}
\def\bp{\begin{proposition}}
\def\ep{\end{proposition}}
\def\be{\begin{equation}}
\def\ee{\end{equation}}

\def\d{{\rm d}}
\def\bea{\begin{eqnarray}}
\def\eea{\end{eqnarray}}
\def\beas{\begin{eqnarray*}}
\def\eeas{\end{eqnarray*}}
\def \pa{\partial}

\def\C{{\mathbb C}}

\def\R{{\mathbb R}}
\def\N{{\mathbb N}}

\date{}
\def \Supp{S_w}
\begin{document}
%

\baselineskip 16pt plus 1pt minus 1pt
\begin{flushright}
CRM-???? (2009)\\
\end{flushright}
\vspace{0.2cm}
\begin{center}
\begin{Large}
\textbf{On the norms and roots of orthogonal polynomials in the plane and $L^p$-optimal polynomials with respect to varying weights}
\end{Large}\\
\bigskip
\bigskip
\begin{large}
{F. Balogh}$^{\dagger\ddagger}$\footnote{balogh@crm.umontreal.ca},
{M. Bertola}$^{\dagger\ddagger}$\footnote{Work supported in part by the Natural
    Sciences and Engineering Research Council of Canada (NSERC).}\footnote{bertola@crm.umontreal.ca}
\end{large}

\bigskip
\begin{small}
$^{\dagger}$ {\em Centre de recherches math\'ematiques,
Universit\'e de Montr\'eal\\ C.~P.~6128, succ. centre ville, Montr\'eal,
Qu\'ebec, Canada H3C 3J7} \\
\smallskip
$^{\ddagger}$ {\em Department of Mathematics and
Statistics, Concordia University\\ 1455 de Maisonneuve Blvd. West, Montr\'eal, Qu\'ebec,
Canada H3G 1M8} \\
\end{small}
\end{center}
\bigskip

\hrule
\begin{center}{\bf Abstract}\\
\end{center}
For a measure on a subset of the complex plane we  consider $L^p$-optimal weighted polynomials, namely, monic polynomials of degree $n$ with a varying weight of the form $w^n = {\rm e}^{-n V}$ which minimize the $L^p$-norms, $1 \leq p \leq \infty$. It is shown that eventually all but a uniformly bounded number of the roots of the $L^p$-optimal polynomials lie within a small neighborhood of the support of a certain equilibrium measure; asymptotics for the $n$th roots of the $L^p$ norms are also  provided. The case $p=\infty$ is well known and corresponds to weighted Chebyshev polynomials; the case $p=2$ corresponding to orthogonal polynomials as well as any other $1\leq p <\infty$ is our contribution.
\medskip \hrule
\bigskip

\section{Introduction, background and results}

In approximation theory an important role is played by the so-called \emph{Chebyshev polynomials} associated to a compact set $K \subseteq \C$, namely monic polynomials of degree $n$ that minimize the supremum norm over $K$. As a natural generalization, one can consider \emph{weighted Chebyshev polynomials} with respect to a varying weight of the form $w^n$ on some $\Sigma \subseteq \C$ that are minimizing the supremum norm of weighted polynomials $Q_n w^n$ over $\Sigma$, where $Q_n$ is a monic polynomial of degree $n$ (the weight function $w$ is assumed to satisfy certain standard admissibility conditions that make the extremal problem well-posed \cite{SaffTotik}).

Along the same lines, given a positive Borel measure $\sigma$ on $\Sigma \subseteq \C$, one can consider optimal weighted polynomials in the $L^2(\sigma)$-sense; provided that the integrals below are finite, it is easy to see that there is a unique monic polynomial $P_n$ for which the weighted polynomial $P_n w^n$ minimizes the $L^2(\sigma)$-norm
\be
\le\|P_n w^n\ri\|_{L^2(\sigma)} := \le(\int_{\Sigma} |P_n|^2 w^{2n} \d \sigma \ri)^{\frac 1 2 }
\ee
among all monic  weighted polynomials of degree $n$. This polynomial may be characterized as the $n$th monic \emph{orthogonal polynomial} with respect to the varying measure $w^{2n} d\sigma$, satisfying
\be
\int_{\Sigma} P_n(z) \ov{z}^k w^{2n}(z)d\sigma(z) = \delta_{kn}h_n \qquad 0 \leq k \leq n\ ,
\ee
where
\be
h_n = \inf\left\{ \le\|Q_n w^n\ri\|_{L^2(\sigma)}\ \colon\ Q_n(z) \mbox{ monic polynomial of degree }n \right\}\ .
\ee

Orthogonal polynomial sequences for varying measures of the form $w^n\d\sigma$ appear naturally in the context of random matrix models \cite{Mehta, Deift}: on the space of $n \times n$ Hermitian matrices ${\mathbb H}_n$, probability distributions of the form
\be
\label{matrint}
\rho_{n}(M)\d M =\frac{1}{{\mathcal Z}_{n}}\exp(-n\tr(V(M)))\d M\ , \qquad  {\mathcal Z}_{n} =\int_{{\mathbb H}_n} \exp(-n\tr(V(M)))\d M
\ee
are considered where the potential function $V(x)$ grows sufficiently fast as $|x| \to \infty$ to make the integral in \ref{matrint} finite ($\d M$ stands for the Lebesgue measure on ${\mathbb H}_n$). The apparent unitary invariance of \ref{matrint} implies that the analysis of statistical observables of $M$ may be reduced to that of the random eigenvalues $\lambda_1,\dots,\lambda_n$ with probability distribution
\be
\label{eigenint}
\begin{split}
&p_{n}(\lambda_1,\dots,\lambda_n) = \frac{1}{Z_n} \prod_{1 \leq k < l \leq n}(\lambda_k-\lambda_l)^2e^{-n\sum_{k=1}^{n}V(\lambda_k)}\ ,\\
&\quad Z_{n} = \int\!\!\cdots\!\!\int_{\R^n}\prod_{1 \leq k < l \leq n}(\lambda_k-\lambda_l)^2e^{-n\sum_{k=1}^{n}V(\lambda_k)}\d\lambda_1\cdots\d\lambda_n\ .
\end{split}
\ee
The marginal distributions of $p_n$ (referred to as \emph{correlation functions}) are expressible as determinants of the weighted polynomials
$p_n(x)e^{-nV(x)/2}$ where the $p_n$ satisfies the orthogonality relation
\be
\int_{\R}p_n(x)x^ke^{-nV(x)}dx = \delta_{kn}h_n\qquad k=0,\dots,n\ .
\ee
Therefore the asymptotic analysis of the correlation functions reduces to the study of the corresponding orthogonal polynomials. On the real line, the asymptotic analysis is done effectively by the so-called Riemann--Hilbert method \cite{Deift}; however, for the so--called {\em normal matrix models} \cite{Zabrodin, ElbauFelder}, for which the eigenvalues may fill regions of the complex plane, much less is known in general. While random  matrix theory was the original impetus behind our interest, the paper will not draw any conclusions on these important connections.

Following instead a more approximation-theoretical spirit, it is also natural to consider {\bf $L^p$-optimal} weighted polynomials \cite{Widom69, LubinskySaff, MhaskarSaff_zeros, MhaskarSaff_lp} with respect to the varying weight $w^n$ and the measure $\sigma$, i.e. to minimize the $L^p$-norm
\be
\le\|P_n w^n\ri\|_{L^p(\sigma)} := \le(\int |P_n|^p w^{np} \d \sigma \ri)^{\frac 1 p }
\ee
over all monic polynomials of degree $n$. 
The paper addresses two questions; the first concerns the location of the roots of $L^p$-optimal polynomials or rather where the roots cannot be. We find that eventually (i.e. for sufficiently large $n$) all roots fall in an arbitrary neighborhood of the convex hull of the support $\Supp$ of the relevant equilibrium measure $\mu_w$ (whose definition is recalled in Sect. \ref{potential}); this is accomplished in Prop. \ref{propconvex} (with a more precise statement).

If the support is not convex (possibly with holes and several disjoint connected components) then we can state that (Prop. \ref{propRunge}) all but a finite (and uniformly bounded) number of roots falls within any arbitrary neighborhood of the polynomially convex hull of the support.
A consequence of the above is that 
\be
\mbox{\Large ``\ }\lim_{n\to\infty} \frac 1n \ln P_n(z) =  \int \ln |z-t|\d\mu_w(t)\mbox{\ \Large ''}
\ee
where the quotation marks indicate that the statement is imprecise (see Thm. \ref{main} for the precise one); the convergence is uniform over closed subsets of the unbounded component  $\C \setminus \Supp$. If $K$ does not contain roots of $P_n$ (eventually) then we can remove the quotations and the statement is correct (for example, if $K$ is disjoint from the convex hull of $\Supp$).

The second question deals with the leading order behaviour of the $L^p$ norms of the $p$--optimal polynomials and we show that -- in fact -- they all have the exact same asymptotic behaviour
\be
\lim_{n\to\infty} \left(\le\|P_n w^{n}\ri\|_{L^{p}(\sigma)}\right)^{1/n} = \exp(-F_w)
\ee
where $F_w$ is the \emph{modified Robin's constant} of the equilibrium measure $\mu_w$, and this limit is independent of $1\leq p\leq \infty$. 

The case of $p=\infty$ of the above statements  is known in the literature (\cite{Fejer} for the unweighted case, and \cite{SaffTotik} for the weighted one)  even in the varying weight case. It seems to be new for $p\neq \infty$.

\subsection{Potential-theoretic background}
\label{potential}
We will consider  polynomials  on a closed set $\Sigma \subseteq \C$, called a {\bf condenser}; on this set a reference measure $\sigma$ is supposed to be given. Since we are not seeking the greatest generality (at cost of simplicity) we will restrict ourselves to the following situations:
\begin{itemize}
\item $\Sigma$ is a finite collection of Jordan  curves, with typical example the real axis or union of intervals thereof. In this case the measure $\sigma$ is simply the arc-length,
\item $\Sigma$ is a finite union of regions of the plane, with the area measure $d\sigma = dA$,
\item $\Sigma$ is a finite union of  elements of both types above.
\end{itemize}
\color[rgb]{0,0,0}
The main focus will be $\Sigma=\C$ or $\Sigma =\R$ or $\Sigma=\gamma$ a smooth curve in $\C$. 

The {\bf weight function} $w\ \colon\ \Sigma \to [0,\infty)$ introduced above is assumed to satisfy the following standard {\bf admissibility} conditions (\cite{SaffTotik}):
\begin{itemize}
\item $w$ is upper semi-continuous,
\item $\mbox{cap}\left(\{z\ \colon \ w(z) > 0\}\right)$ has positive capacity,
\item $|z|w(z) \to 0$ as $|z|\to \infty$ in $\Sigma$.
\end{itemize}
The {\bf potential} $V(z)$ is the function for which $w(z) = \exp(-V(z))$ and it inherits the corresponding admissibility conditions.
The weighted energy functional is defined as follows; for a probability measure $\mu$ on $\Sigma$ we define
\be
\mathcal I_w[\mu]:=  \int\!\!\int \ln \frac 1{|z-w|} \d\mu(z)\d\mu(w)+ 2\int V(z)\d\mu(z)\ .\label{Energy}
\ee
It is well known in potential theory \cite{SaffTotik} that there exists a unique measure $\mu_w$ that realizes the minimum of ${\mathcal I}_w$; such a  measure is referred to as the {\bf equilibrium measure}. Its support $\Supp={\rm supp}(\mu_w)$ is a compact set.

Although it will not be used directly we recall the following indirect characterization of $\mu_w$: if we denote with 
\be
U^{\mu}(z):= \int \ln \frac 1{|z-w|}\d\mu(w)
\ee
the {\bf logarithmic potential} of a probability measure $\mu$ then $\mu_w$ is uniquely characterized as follows. There exists a constant $F_w$ called the {\bf modified Robin's constant} such that the {\bf effective potential}
\be
\Phi(z):= U^{\mu_w}(z)+V(z) -F_w  
\ee
satisfies 
\be
\left\{
\begin{split}
\ds \Phi(z)\leq 0& \quad z\in \Supp\\
& \mbox{ and }\\
\ds \Phi(z)\geq 0& \quad  z\in \Sigma \ \quad \mbox{q. e.}
\end{split}
\right\}
 \Rightarrow \Phi(z)= 0  \quad z\in \Supp\quad \mbox{q. e.}
\ee
where 'q. e.' stands for ``quasi-everywhere'', namely up to sets of zero logarithmic capacity.

\section{Where the roots are not}
Let $P_n(z)$ be  any sequence of  polynomials of degree $\leq n$,  $\Supp=\supp(\mu_w)$ and let $\mathcal N\supset \Supp$ be an open bounded  set containing $\Supp$. 

In \cite{SaffTotik} III.6 (eq. 6.4) it is shown in general (under certain assumptions on $\Sigma, w$ and $\sigma$) that if $P_n$ is any sequence of polynomials of degree $\leq n$ we have
\be
\|P_n w^n\|_{L^p(\sigma)}^p = \int_{\Sigma} |P_nw^n|^p \d\sigma \leq (1+ C {\rm e}^{-c n })\int_{\mathcal N} |P_nw^n|^p \d\sigma \label{311}
\ee
where the constants $c>0$ and $C$ {\bf do not depend} on the polynomial sequence under consideration (they depend --of course-- on $w, p$ and  $\mathcal N$).

The inequality (\ref{311}) can be rewritten 
or equivalently  ($\chi_{_\mathcal N}$ denotes the indicator function of the set $\mathcal N$)
\be
1 \leq \frac {\|P_n w^n\|_p^p}{\|P_nw^n\chi_{_{\mathcal N}}\|_p^p}\leq 1+ C{\rm e}^{-cn}\ .
\ee

\begin{wrapfigure}{r}{0.5\textwidth}
 \resizebox{0.4\textwidth}{!}{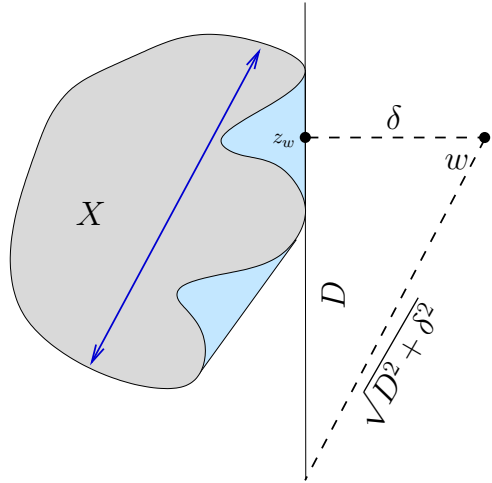}
 \caption{Figure for Lemma \ref{lemmaconvex}}
 \label{convex}
 \end{wrapfigure}

The inequality (\ref{311}) shows that the norm of $P_nw^n$ lives in a small neighborhood of $\Supp$; this will be the main tool in what follows.
The ideas follow very closely similar steps for the so--called {\em weighted Chebyshev polynomials} in III.3 of \cite{SaffTotik}.

For any set $X\subset \C$ we will denote by  $\wh X$ the {\bf (closed) convex hull} of said set.

 Let, as before $\mathcal N\supset \Supp$ be  an open, bounded neighborhood of $\Supp$. We start from the 

 \bl
 \label{lemmaconvex}
 Let  $X\subset \C$ be compact that is not a singleton and $w\in \C$ be such that ${\rm dist}(w,\wh X )=\delta>0$. Then 
 \be
 \frac {|z-z_w|}{|z-w|}\leq \frac {D}{\sqrt{D^2  + \delta^2}}<1 \ ,\qquad D:= {\rm diam}( \wh X)
 \label{26}
 \ee
 where $z_w\in \wh{X}$ is  the (unique) closest point to $w$.
 \el
 {\bf Proof.} 
 The set $\wh {X}$ lies entirely on one half-plane passing through $z_w$ and perpendicular to the line segment $[z_w,w]$.
 Let $\theta_w$ the smallest angle such that $\wh{X}$ is entirely contained in a $\theta_w$ sector centered at $w$; by the convexity and compactness of $\wh{X}$, $\theta_w<\pi$. In fact we can estimate the upper bound on $w$ of such $\theta_w$ as 
 \be
 \theta_w \leq \arctan\le(\frac {\delta}{D}\ri)\ ,\qquad D = {\rm diam} (\wh X)\ .
 \ee
 from which (\ref{26}) follows (see Fig. \ref{convex}).
 \QED

\bp
\label{propconvex}
Let $K$ be a closed subset  in $\C\setminus \wh \Supp$. Then eventually there are no roots of $P_n$ belonging to  $K$. In particular, for any $\epsilon>0$ there is a $n_0\in\N$ such that $\forall n>n_0$ all roots of $P_n$ are within distance $\epsilon$ from the convex hull.
\ep
{\bf Proof.}
Since $K$ is closed and has no intersection with $\wh \Supp$ we have ${\rm dist}(K,\wh \Supp) = 2\delta>0$;
Let $\mathcal N$ be the $\delta$-fattening of $\wh \Supp$, namely
\be
\mathcal N := \{z\in \C\ \colon \ {\rm dist}(z,\wh \Supp)\leq  \delta\}
\ee
It is easy to see that $\mathcal N$ is convex as well. 

Now consider the $p$--optimal polynomial $P_n(z)$ and let us decompose it as $P_n(z) = R_n(z) Q_n(z)$ where $R_n(z)$ is the factor of all roots within $K$; note that each of these roots is at distance $\geq \delta$ from $\mathcal N$. 

For each root $z_j$ of $R_n(z)$ we can find the closest point $\wt z_j\in \mathcal N$; hence we will define $\wt R_n(z)$ as the ``proximal substitute'' of $R_n$, where each root of $R_n$ has been replaced by its proximal point in $\mathcal N$.
Then for all $z\in {\mathcal N}$ we have $|\wt R_n(z)|\leq \rho^{r_n}  |R_n(z)|$ where $r_n= \deg (R_n)$. Indeed, by Lemma \ref{lemmaconvex},
\bea
|\wt R_n(z)| = \prod_{j=1}^{r_n} |z - \wt z_j| \leq\rho^{r_n} \prod_{j=1}^{r_n} |z-z_j| =\rho^{r_n} |R_n(z)| 
\eea
Thus {\bf pointwise}
\be
|\wt P_n(z)| \leq  \rho^{r_n} |P_n(z)|\ \ z\in\mathcal N\ ,\ \ \rho:= \frac {D}{\sqrt{D^2+\delta^2}}<1 \ ,\ \ D:={\rm diam}(\mathcal N).
\label{212}
\ee

We thus have
\bea
\le \|\wt P_nw^n\chi_{_{\mathcal N}} \ri \|_{p}^p 
\mathop{\leq}^{\mbox{\tiny{by (\ref{212})}}} \rho^{pr_n}\le \|P_nw^n\chi_{_{\mathcal N}} \ri\|_{p}^p\leq \rho^{pr_n}  \le\|P_n w^n\ri \|_p^p \label{goodone}
\eea
By definition, the $p$--optimal polynomials $P_n$ have the smallest $L^p$ norm and hence 
\be
1\leq \frac { \le\| \wt P_n w^n  \ri\|_p^p}{ \le\|P_n  w^n \ri\|_p^p} 
\mathop{\leq}^{\mbox{(\ref{311})}}   (1+ C{\rm e}^{-cn})\frac { \le\| \wt P_n w^n \chi_{_{\mathcal N}} \ri\|_p^p}{ \le\|P_n  w^n \ri\|_p^p} \mathop{\leq}^{\mbox{(\ref{212})}} (1+C{\rm e}^{-cn}) \rho^{pr_n}.\label{verygood}
\ee
where in the second inequality we have used (\ref{311}) on the sequence of polynomials $\wt P_n$. Inequalities (\ref{verygood}) amount to 
\be
1\leq (1+C{\rm e}^{-cn}) \rho^{p r_n}
\ee
and recall that $\rho<1$. This inequality implies at once that $\limsup r_n = 0$, and hence the sequence of natural numbers $r_n$ must eventually be identically zero. 
The second statement in the theorem is simply obtained by taking for $K$ the complement of the $\epsilon$--fattening of ${\rm Co}(\Supp)$.
\QED

Having established that there are no roots (eventually) ``outside'' the convex hull, we get some further information about what happens in general.

 We borrow the following nice 
\bl[Lemma III.3.5 in \cite{SaffTotik}, originally in \cite{Widom67}]
\label{lemmaRunge}
If $S$ and $K$ are compact sets such that $\mbox{Pc}(S) \cap K =\emptyset$ then there is a positive integer $m=m(K)$ and a constant $0<\alpha(K)<1$ such that for all $(z_1,\dots,z_m)\in K^m$ there are points $\wt z_1, \dots, \wt z_m$ such that the rational function
\be
r(z):= \frac{\prod_{j=1}^{m(K)}  (z-\wt z_j)}{\prod_{j=1}^{m(K)} (z-z_j)} 
\ee
satisfies 
\be
\sup_{z\in S} |r(z)|\leq \alpha(K)\ .
\ee
\el

Lemma \ref{lemmaRunge} allows us to prove 
\bp
\label{propRunge}
For any compact set $K$ contained in the unbounded component of $\C\setminus \Supp$ the number of roots of the $p$--optimal  polynomials $P_n$ within $K$ is bounded. In particular $\forall \epsilon>0$ there is $n_0\in \N$ such that $\forall n>n_0$ all but a finite number (uniformly bounded) roots of $P_n$ lie within distance $\epsilon$ from the polynomial convex hull of $\Supp$ (i.e. $\C \setminus \Omega$, where $\Omega$ is the unbounded component of $\C \setminus \Supp$).
\ep
{\bf Proof.}
In parallel with the proof of Prop. \ref{propconvex} let $2\delta = {\rm dist}(K,\Supp)$ and let $\mathcal N$ be the  $\delta$-fattening of $\Supp$. We decompose $P_n= R_n Q_n$ where $R_n$ has $r_n$ roots (counted with multiplicity) within $K$. We will prove that  $r_n < m(K)$ eventually, where $m(K)$
is the number of poles in Lemma \ref{lemmaRunge} for $S=\mathcal N$ and $K$.
Proceeding by contradiction, there would be a subsequence where $r_n \geq m(K)$; but then we can use Lemma \ref{lemmaRunge} to find  a polynomial $\wt R_n$ such that 
\be
|\wt R_n(z)|\leq \alpha(K) |R_n(z)|\ ,\qquad z\in \mathcal N\ \Rightarrow \ \  |\wt P_n(z)| \leq \alpha(K) |P_n(z)| \ , \ \ z\in \mathcal N.
\label{216}
\ee
At this point we proceed exactly as in the proof of Prop. \ref{propconvex} starting from (\ref{goodone}) with $\rho^{p r_n} \mapsto \alpha(K)$, namely,
\bea
\le \|\wt P_nw^n\chi_{_{\mathcal N}} \ri \|_{p}^p 
\mathop{\leq}^{\mbox{\tiny{by (\ref{216})}}} \alpha(K)^p \le \|P_nw^n \chi_{_{\mathcal N}} \ri\|_{p}^p\leq \alpha(K)^p \le\|P_n w^n\ri \|_p^p \label{goodone2}
\eea
By the $p$--optimality of the polynomial  $P_n$ we must have 
\be
1\leq \frac { \le\| \wt P_n w^n  \ri\|_p^p}{ \le\|P_n  w^n \ri\|_p^p} 
\mathop{\leq}^{\mbox{(\ref{311})}}   (1+ C{\rm e}^{-cn})\frac { \le\| \wt P_n w^n \chi_{_{\mathcal N}} \ri\|_p^p}{ \le\|P_n  w^n \ri\|_p^p} \mathop{\leq}^{\mbox{(\ref{216})}} (1+C{\rm e}^{-cn}) \alpha(K)^p.\label{verygood2}
\ee
It is clear that the last expression in (\ref{verygood2}) is eventually {\em less than one} (since $\alpha(K)<1$), which leads to a contradiction with the assumption that there were $\geq m(K)$ roots in $K$.
The last statement follows from the fact that there are no roots outside the convex hull by Prop. \ref{propconvex} together with the above.
 \QED
 \bx
 Suppose that the support of the equilibrium measure consists of intervals in the real axis, as in the case of ordinary orthogonal polynomials.
 It is an exercise to see that for any gap the number $m(K)=2$ and hence there can be at most one zero within each gap.
 \ex
We next prove

\bt
\label{main}
Let $\Omega$ be the unbounded connected component of $\C \setminus \Supp$ and $K\subset \Omega$ a compact subset.
Let $z_{\ell,n}(K)$ be the roots of $P_n$ belonging to $K$,  $\ell = 1,\dots, m_n(K)$.
Then, uniformly in $K$ we have 
\be
\lim_{n\to\infty} \frac 1 n \ln |P_n(z)| + \frac 1 n \sum_{j=1}^{m_n(K)} G_{\Omega}(z,z_{\ell,n}) = \int \ln |z-t|\d\mu(t)\ ,
\ee
where $G_\Omega(z,w)$ is the Green's function of $\Omega$, namely the function such that 
\bea
&& \triangle_z G_{_\Omega}(z,w) \equiv 0 \ ,\ \ z\in \Omega \setminus \{w\}\\
&& G_{_{\Omega}}(z,w) =0\ ,\ \ z\in \pa\Omega\\
&& G(z,w)\geq 0\ ,\ \ z,w\in\Omega\\
&& G_{_{\Omega}}(z,w) = \ln \frac 1 {|z-w|} + \mathcal O(1)\ \ z\to w\label{Gasym}
\eea
Additionally, if $K$ is closed and does not contain (eventually) any roots, then, uniformly, 
\be
\lim_{n\to\infty} \frac 1 n \ln |P_n(z)|  = \int \ln |z-t|\d\mu(t)\label{Kasym}
\ee
\et
{\bf Proof.}
We reason on the functions 
\be
f_n(z):=  \frac 1 n \ln |P_n(z)| + \frac 1 n \sum_{j=1}^{m_n(K)} G_{\Omega}(z,z_{\ell,n}) -  \int \ln |z-t|\d\mu(t)
\ee
We will see in Prop. \ref{lowerest} together with Corollary \ref{coras} that  $\forall \epsilon>0$ $\exists n_0: \ n\geq n_0$
\bea
\frac 1 n \ln |P_n(z)w^n(z)| \leq -F_w + \epsilon\ ,\qquad \forall z\in \C
\eea
Additionally, the $f_n(z)$'s  are subharmonic in $\Omega$ and harmonic in a neighborhood of  $z=\infty$: indeed all roots are uniformly bounded (from Prop. \ref{propconvex}) and  the Green's function $G_\Omega(z,w)$ is harmonic away from the singularity $z=w$ (in a neighborhood of which it is {\em super}harmonic) and in the neighborhoods of $z_{\ell,n}$ the $f_n$'s are actually harmonic because the singularities coming from $P_n$'s cancel out exactly those coming from the Green's functions.

For  $z\in\pa \Omega$ and  $\forall \epsilon>0$ we have eventually (recall that $G_\Omega(z,w)=0$ for $z\in\pa \Omega$)
\be
f_n(z) \leq V(z) +U^{\mu_w}(z) -F_w + \epsilon  \leq \epsilon\ ,\qquad z\in \pa\Omega.
\ee
Since $f_n(z)$ are subharmonic, they cannot have isolated maxima in the interior of $\Omega$ and hence we conclude that $f_n(z)\leq \epsilon$ throughout $\Omega$ (including $z=\infty$). 

 Let $f_\infty(z)=\limsup_{n\to\infty} f_n(z)$; then $\forall \epsilon>0$ 
\be
f_\infty(z) = \limsup_{n\to\infty} f_n(z) \leq \epsilon\ \Rightarrow f_\infty(z)\leq 0 \ ,\qquad z\in \C\label{done}
\ee
Let $\mbox{Pc}(\Supp) = \C\setminus \Omega$ be the {\em polynomial convex hull} of $\Supp$ and 
let now $K$ be a {\bf compact} set $K\subset \Omega$.

We next analyze the $\liminf$; let $z_0\in K$ and set 
\be
L_{z_0}:= \liminf f_n(z_0) \leq f_\infty(z_0)\leq 0\ .
\ee 
where $z_0\in K $ is some (arbitrary but fixed) point.
There is a subsequence $n_k$ of the numbers $f_n(z_0)$'s which converges to this limit; out of it, we can extract another subsequence (which we denote again $n_k$ for brevity) such that the counting measures $\sigma_{n_k}$ have a weak$^\star$ limit (since they are all compactly supported) which we denote by $\sigma_{z_0}$ (note that both the subsequence and this limiting distribution may depend on $z_0$). Prop. \ref{propRunge} implies that its support of $\sigma_{z_0}$ lies in the polynomial convex hull of $\Supp$; in particular the function $\ln |z-\bullet |$ is harmonic on ${\rm supp}(\sigma_{z_0})$ for any $z\in K$.   Let $\widehat \sigma_{n_k}$ be the restriction of $\sigma_{n_k}$ to those atoms outside of $K$; we know that it differs from $\sigma_{n_k}$ by a finite number $m_n(K)$ of atoms (uniformly bounded in $n$) and hence it obviously has the same weak$^\star$ limit.  Now, for any $z\in K$ along the chosen subsequence we have 
\bea
0\geq f_\infty(z)\geq \lim_{k\to\infty} f_{n_k}(z)= \lim_{k\to\infty} \int \ln |z-t|\d\widehat\sigma_{n_k} (t) - \int \ln|z-t|\d\mu_w(t) +\cr + \frac 1 {n_k}  \sum_{\ell=1}^{m_n(K)} \le(G_{_\Omega}(z,z_{\ell,n_k}) + \ln|z-z_{\ell,n_k}|\ri) \label{224}
\eea
Since $G_{_\Omega}(z,w)+ \ln|z-w|$ is jointly continuous in $z,w$ for $z,w\in \Omega$, it is also (jointly) bounded on compact sets; we know already that $z_{\ell,n_k}$ all are uniformly bounded, hence the last term in (\ref{224}) tends to zero. We thus have
\bea
\lim_{k\to\infty} f_{n_k}(z)&\&= \lim_{k\to\infty} \int \ln |z-t|\d\widehat\sigma_{n_k} (t) - \int \ln|z-t|\d\mu_w(t) = \cr
&\& =  \int \ln |z-t|\d\widehat\sigma_{z_0} (t) - \int \ln|z-t|\d\mu_w(t) 
\label{224bis}
\eea
The right hand side of (\ref{224bis}) is harmonic in $\Omega$ (by inspection)  and by (\ref{done}) it is $\leq 0$; on the other hand at $z=\infty$ it vanishes (since both measures are probability measures) and hence it must be identically zero. Evaluating it at $z=z_0$ yields that $L_{z_0} = \liminf_{n\to\infty} f_n(z_0)=0$; since $z_0$ was arbitrary, this shows that $\lim_{n\to\infty} f_n(z)=0$; the uniformity of the convergence follows from the fact that the sequence of functions 
\be
h_n(z):= \int\ln|z-t| \d \widehat \sigma_n(t)
\ee
are {\em equicontinuous} for $z\in K$ and hence the Arzela--Ascoli Theorem \cite{Rudin} guarantees uniform convergence. To see equicontinuity we compute
\bea
\le|h_n(z)-h_n(z') \ri| = \int \ln \le|\frac {z-t}{z'-t}\ri| \d \widehat \sigma_{n_k}(t)  =  \int \ln \le| 1 + \frac {z-z'}{z'-t}\ri| \d \widehat \sigma_{n_k}(t)  \leq  \int  \frac {|z-z'|}{|z'-t|} \d \widehat \sigma_{n_k}(t)\leq \cr
\leq \frac {|z-z'|}{{\rm dist}(K,\Supp)}
\eea
Note that the above chain of inequalities applies more generally for any {\em closed} $K\subset \Omega$ and also says that the sequence is uniformly Lipschitz.

To prove (\ref{Kasym}) we note that we have used compactness only after (\ref{224}), but if $m_n(K)\equiv 0$ (eventually) then the same arguments prove uniform convergence without having to use compactness.
\QED

Theorem \ref{main} says loosely speaking that $\frac 1n \ln |P_n|$ converges to the logarithmic transform of the equilibrium measure as uniformly as it is possible on the ``outside'' of the support, given that there are possibly some stray roots; if we restrict to the outside of the convex hull of $\Supp$, then this convergence is truly uniform (over closed subsets) because --eventually--  there are no roots at all (Prop. \ref{propconvex}).

Theorem \ref{main} has an interesting corollary 
\bc
\label{corbalayage}
Let $h(z)$ be any harmonic function on a neighborhood of $\mbox{Pc}(\Supp)$ and let $\sigma$ be a weak$^*$ limit point of the counting measures of the $L^p$--optimal polynomials. Then 
\be
\int h(z) \d\sigma(z) = \int h(z)\d\mu_w(z)\ .
\ee
\ec
{\bf Proof.}
By Mergelyan's theorem it suffices to verify it for the monomials $z^j$; we have seen in the proof of Thm. \ref{main} (\ref{224} and discussion thereafter) that 
\be
\int \ln |z-t|\d\sigma(t) - \int \ln|z-t|\d\mu_w(t) \equiv 0
\ee
for $z\in \Omega$ (the complement of the polynomial convex hull of $\Supp$). Taking the large $z$ expansion we have easily the statement
\QED
\br
The Theorem \ref{main} and Corollary \ref{corbalayage} assert that whatever limiting distribution the roots of the $p$--optimal polynomials may have, it must be  a balayage of the equilibrium measure onto the support of this limiting distribution. In order not to swindle the reader, we should point out that it falls short of saying that there is a {\em unique} limiting distribution, and even further away from any statement about what distribution that should be.
\er


\section{Norm estimates}

\subsection{Upper estimate for the norms}
The aim of this section is twofold:
first we will prove that if $P_nw^n$ are the $p$-optimal weighted polynomials then 
\be
\lim_{n\to\infty} \frac 1 n \ln\|P_nw^n\|_p = -F_w \ \ \Leftrightarrow\ \  \|P_nw^n\|_p = {\rm e}^{-nF_w + o(n)}
\ee
{\em En route} we will see that the $L^p$ norms of the wave-functions $P_nw^n$ are asymptotically equal to the $L^\infty$ ones. In particular this implies that the $n$-th root of the wave functions is uniformly bounded.
\bp
\label{upperest}
Let $P_nw^n$ be the $p$--optimal weighted polynomial; then 
\be
\limsup_{n\to \infty} \frac 1n \ln \|P_nw^n \|_p \leq -F_w\ ,
\ee
where $\ell$ is the Robin constant for the equilibrium measure.
\ep 
{\bf Proof.}
We compare the $L^p$ norms of the $P_nw^n$'s  with the weighted Fekete polynomials $F_n w^n$.  Let $\mathcal N $ be a bounded open neighborhood of $\Supp$.
Then 
\bea
\|P_nw^n\|_p \mathop{\leq}^{\mbox{\tiny by optimality}} \|F_nw^n \|_p \mathop{\leq}^{\mbox{(\ref{311})}} 
 (1+ C{\rm e}^{-cn}) \|F_nw^n \chi_{\mathcal N}\|_p\leq   (1+ C{\rm e}^{-cn})\| F_n w^n\|_\infty  {\mathcal Area}(\mathcal N)^{\frac 1 p} 
\eea
Now taking $\frac 1n \ln(\cdot)$ of both sides gives
\be
\frac 1 n \ln \le(\|P_nw^n\|_p\ri)\leq \frac 1 n \ln \le(\|F_n w^n\|_p\ri) 
\leq \frac 1n \ln \le(\| F_n w^n\|_\infty\ri)  + \mathcal O(n^{-1})
\ee
Since 
\be
\lim_{n\to \infty}  \frac 1n \ln \le(\| F_nw^n\|_\infty\ri) = -F_w
\ee
(see Thm. III.1.9 in \cite{SaffTotik})  we have 
\be
\limsup_{n\to \infty} \frac 1 n \ln \le(\|P _nw^n \|_p\ri)\leq \limsup_{n\to \infty} \frac 1 n \ln \le(\|F_nw^n \|_p\ri) \leq \limsup_{n\to \infty} \frac 1 n \ln \le(\|F_nw^n \|_\infty\ri) \leq -F_w\ .
\ee 
\QED
\br
It may be of some importance to note that the above proof can be used to show\ .
\be
\limsup_{n \to \infty} \ln \frac{\|P_n w^n\|_p}{\|F_nw^n\|_\infty}\leq \sqrt[p]{{\mathcal Area}(\Supp)}
\ee
\er
\subsection{Lower estimate for the norms}
We follow the idea in \cite{SaffTotik}, pp 182.

\bl
\label{lowerest}
Let $P_n(z)$ be a sequence of polynomials of degree at most $n$. Assume further that the potential $V$ is twice continuously differentiable. Then there is a constant $D>0$ and $d_\Sigma$ (the Hausdorff dimension of $\Sigma$, which for us is either $2$ or $1$) such that 
\be
\frac { \|P_n w^n \|_p}{\|P_nw^n\|_\infty}  \geq D  n^{-\frac{d_\Sigma}p}
\ee
In particular 
\bea
\liminf_{n\to \infty} {\|P_n w^n \|_{p}}^{\frac 1n}\geq \liminf_{n\to\infty} { \|P_n w^n \|_{\infty}}^{\frac 1n}\label{439}\\
\limsup_{n\to \infty} {\|P_n w^n \|_{p}}^{\frac 1n}\geq \limsup_{n\to\infty} { \|P_n w^n \|_{\infty}}^{\frac 1n}\ .
\eea
\el
{\bf Proof.}
We work with the normalized polynomials
\be
Q_n(z):= \frac 1{ \|P_nw^n\|_\infty} P_n(z).
\ee

Let $z_0$ be a point where $|Q_n(z)w^n(z)|$ achieves its maximum value $1$ (such a point exists by the assumed admissibility conditions on $w$). 
We claim that 
\bea
\exists C>0: \ |z-z_0|\leq \frac 1{2{\rm e} C n} \  \Longrightarrow\ \ |Q_n(z)| {\rm e}^{- n V(z)} \geq \frac 1 {2{\rm e}}
\label{smalldisk}
\eea
Since $|Q(z_0)|{\rm e}^{-nV(z_0)}=1$ the inequality can be rewritten
\bea
|Q_n(z)|\leq|Q_n(z_0)|  {\rm e}^{n (V(z)-V(z_0))} \ ,\ \ \forall z\in \C.
\eea
Let $\delta>0$ and set $C_\delta(z_0):= \sup_{|z-z_0|=\delta} |V(z)-V(z_0)|$; since we are assuming $V(z)$ to be twice continuously differentiable,  $z_0\in \Supp$ and $\Supp$ is compact,   we see that a simple argument shows  $C_\delta(z_0)< C \delta$ for some constant $C>0$ (independent of $z_0\in \Supp$).
Let $|z-z_0|<\frac 1 2  \delta$; the formula of Cauchy for the derivative implies 
\be
|Q'_n(z)| \leq  |Q_n(z_0)|\frac 2{\delta}{\rm e}^{n C \delta }\ ,\ \ \ |z-z_0|\leq \frac 1 2  \delta.
 \ee
On the even smaller disk $|z-z_0|<\frac1 {4{\rm e}} \delta$ we have 
 \be
 |Q_n(z)-  Q_n(z_0)| \leq \int_{z_0}^{z} |Q_n'(t)| |\d t| \leq  |Q_n(z_0)|\frac{2 {\rm e}^{nC \delta  }|z-z_0|}{\delta} \leq
  \frac 1 2 |Q_n(z_0)| {\rm e}^{nC \delta-1  }
 \ee
 If we choose  $\delta = \frac 1{Cn}$  and hence $|z-z_0|<\frac \delta{4{\rm e}} = \frac 1{4 C n {\rm e}}$ we have 
 \be
 |Q_n(z)-  Q_n(z_0)| \leq\frac 1 2 |Q_n(z_0)| \ \ \Rightarrow \ \ \ 
 |Q_n(z)| \geq \frac 1 2  |Q_n(z_0)|\ .
 \ee
Multiplying both sides
\bea
 |Q_n(z)| {\rm e}^{-n(V(z)-V(z_0))} &\& \geq 
\frac 1 2  |Q_n(z_0)| {\rm e}^{-n(V(z)-V(z_0))}\geq \frac 1 2 |Q_n(z_0)| {\rm e}^{-n C \delta } = \frac {|Q_n(z_0)|}{2{\rm e}} \ \ \Rightarrow\cr
 &\& |Q_n(z)| {\rm e}^{-nV(z)} \geq  \frac 1{2{\rm e}}|Q_n(z_0)|{\rm e}^{-n V(z_0)} =\frac 1{2{\rm e}} \ .
\eea

Integrating the inequality (\ref{smalldisk}) 
\bea
\le(\int_\Sigma |Q_n(z) w^n |^p \d_\Sigma z \ri)^{\frac 1 p}\geq
\le(\int_{|z-z_0|<\frac {\delta}{4{\rm e}}} |Q_n(z) w^n |^p \d_\Sigma z \ri)^{\frac 1 p} \geq\\
\geq 
\frac 1{2{\rm e}}\le[ B_\Sigma\le(\frac 1{4n C{\rm e}}\ri)\ri]^{\frac  1 p}\label{intsmalldisk}
\eea
Here $B_\Sigma(\delta)$ is the $\d\sigma$ volume of the ball of radius $\delta$ centered at $z_0$ in $\Sigma$: in the case $\Sigma=\C$ this is simply $\pi \delta^2$, in the case $\Sigma$ is a smooth curve then $B_\Sigma(\delta) \geq c \delta$ for some $c>0$. The only important fact for us below is that $B_\Sigma(\delta)$ is bounded below by some positive power of $\delta$.
Therefore, recalling that $Q_n (z) = P_n(z)/\|P_nw^n\|_\infty$ the inequality (\ref{intsmalldisk}) reads 
\be 
\le(\|P_n(z) w^n \|_p\ri) \geq \frac 1{2{\rm e}}  
\le[ B_\Sigma\le(\frac 1{4n C{\rm e}}\ri)\ri]^{\frac  1 p}
\|P_nw^n\|_\infty
\ee
Summarizing, there are constants $D>0$ and $d_\Sigma$ (the ``dimension'' of $\Sigma$, which for us is either $2$ or $1$) such that 
\be
\frac { \|P_n(z) w^n \|_p}{\|P_nw^n\|_\infty}  \geq D  n^{-\frac{d_\Sigma}p}\ .
\ee
 \QED
 Before proceeding we recall 
 \bt[Thm. I.3.6 in \cite{SaffTotik}]
 \label{asextr}
 Let $P_n$ be any sequence of {\em monic}  polynomials of degree $n$. Then 
 \be
 \liminf_{n\to \infty} \le(\|  P_n w^n \|_\infty\ri)^{\frac 1n} \geq \exp(-F_w)\ .
\ee
 \et
As a corollary of Thm. \ref{asextr} and Prop. \ref{upperest} we have 
 \bc
 \label{coras}
 The norms of the $p$--optimal polynomials satisfy 
 \be
\sqrt[n]{ \| P_nw^n\|_p }\to {\rm e}^{-F_w}\ , \qquad n \to \infty\ .
 \ee
 \ec
 {\bf Proof.}
Using Lemma \ref{lowerest}  and (\ref{439}) together with 
Thm. \ref{asextr} we have  that the $\liminf$ of the left hand side cannot be less than ${\rm e}^{-F_w}$:
 \beas
 -F_w &\& \mathop{\geq}^{\mbox{\tiny Prop. \ref{upperest}}} \limsup_{n\to\infty} \frac 1 n \ln \|P_n w^n\|_p \geq \liminf_{n\to\infty} \frac 1 n \ln \|P_n w^n\|_p 
  \mathop{\geq}^{\mbox{\tiny Prop. \ref{lowerest}}} \liminf_{n\to\infty} \frac 1 n \ln \|P_n w^n\|_\infty  \mathop{\geq}^{\mbox{\tiny Thm. \ref{asextr}}} -F_w\ .
 \eeas
 \QED
\bibliographystyle{unsrt}
\bibliography{PlaneOProots.bib}

\end{document}

%% file: ConvexN.pdf_t
\begin{picture}(0,0)%
\includegraphics{ConvexN.pdf}%
\end{picture}%
\setlength{\unitlength}{3947sp}%
\begingroup\makeatletter\ifx\SetFigFont\undefined%
\gdef\SetFigFont#1#2#3#4#5{%
  \reset@font\fontsize{#1}{#2pt}%
  \fontfamily{#3}\fontseries{#4}\fontshape{#5}%
  \selectfont}%
\fi\endgroup%
\begin{picture}(4082,4089)(2341,-5617)
\put(4576,-2735){\rotatebox{359.0}{\makebox(0,0)[lb]{\smash{{\SetFigFont{11}{13.2}{\familydefault}{\mddefault}{\updefault}{\color[rgb]{0,0,0}$z_w$}%
}}}}}
\put(6035,-2954){\rotatebox{359.0}{\makebox(0,0)[lb]{\smash{{\SetFigFont{17}{20.4}{\familydefault}{\mddefault}{\updefault}{\color[rgb]{0,0,0}$w$}%
}}}}}
\put(2909,-3424){\makebox(0,0)[lb]{\smash{{\SetFigFont{17}{20.4}{\familydefault}{\mddefault}{\updefault}{\color[rgb]{0,0,0}$X$}%
}}}}
\put(5191,-4138){\rotatebox{90.0}{\makebox(0,0)[lb]{\smash{{\SetFigFont{17}{20.4}{\familydefault}{\mddefault}{\updefault}{\color[rgb]{0,0,0}$D$}%
}}}}}
\put(5446,-5174){\rotatebox{60.0}{\makebox(0,0)[lb]{\smash{{\SetFigFont{17}{20.4}{\familydefault}{\mddefault}{\updefault}{\color[rgb]{0,0,0}$\sqrt{D^2+\delta^2}$}%
}}}}}
\put(5541,-2614){\makebox(0,0)[lb]{\smash{{\SetFigFont{17}{20.4}{\familydefault}{\mddefault}{\updefault}{\color[rgb]{0,0,0}$\delta$}%
}}}}
\end{picture}%